\documentclass[a4paper,11pt,fullpage]{article}
\usepackage[english]{babel}

\usepackage[left=1in,right=1in,top=1in,bottom=1in]{geometry}

\usepackage[usenames]{color}
\usepackage{colortbl}

\usepackage{indentfirst}
\usepackage{amsmath}
\usepackage{amsthm}
\usepackage{amsfonts}
\usepackage{amssymb}
\usepackage{amsxtra}
\usepackage{latexsym}
\usepackage{ifthen}
\usepackage{cite}
\usepackage{esint}
\usepackage[cp1250]{inputenc}

\usepackage{epic}
\usepackage{eepic}

\numberwithin{equation}{section}

\def\XXint#1#2#3{{\setbox0=\hbox{$#1{#2#3}{\int}$}
     \vcenter{\hbox{$#2#3$}}\kern-.5\wd0}}

\begin{document}

\begin{center}{\Large {\bf Maximum principle for the weak solutions of the Cauchy problem for the fourth-order hyperbolic equations}}\\

Kateryna Buryachenko

\end{center}

\begin{abstract}
We investigate the maximum principle for the weak solutions to the Cauchy problem for the hyperbolic fourth-order linear equations with constant complex coefficients in the plane bounded domain. 
\end{abstract}

\section*{Introduction}

We concern here the problem of proving the analog of maximum principle for the weak solutions of the Cauchy problem for the fourth-order linear hyperbolic equations with the complex constant coefficients and homogeneous non-degenerate symbol in some plane bounded domain $\Omega\in R^2$ convex with respect to characteristics:

 \begin{equation}
L(\partial_x)u=a_0\frac{\partial^4u}{\partial x_1^4}+a_1\frac{\partial^4u}{\partial x_1^3\partial x_2}+a_2\frac{\partial^4u}{\partial x_1^2\partial x_2^2}+a_3\frac{\partial^4u}{\partial x_1\partial x_2^3}+a_4\frac{\partial^4u}{\partial x_2^4}=f(x).\label{eq:01}
\end{equation}
 Here coefficients $a_j,\,j=0,\,1,...,\,4$ are constant, $f(x)\in L^2(\Omega),$ $\partial_x=\left(\frac{\partial}{\partial x_1},\,\frac{\partial}{\partial x_2}\right).$ We assume, that Eq. (\ref{eq:01}) is hyperbolic, that means that all roots of characteristics equation
 \begin{equation*}
 L(1,\,\lambda)=a_0\lambda^4+a_1\lambda^3+a_2\lambda^2+a_3\lambda +a_4=0
 \end{equation*}
 are  prime, real  and are not equal to  $\pm i$, that means that the symbol of Eq. (\ref{eq:01}) is non-degenerate or that the Eq. (\ref{eq:01}) is a principal-type equation. The equations for which the roots of the corresponding characteristic equation are multiple and can take the values $\pm i$ are called the equation with
degenerate symbol (see \cite{Bur1}).

 The main novelty of the paper is to prove the analog of maximum principle for the fourth-order hyperbolic equations. This question is very important due to usually a natural physical interpretation, and that it helps to establish the qualitative properties of the solutions (in our case the questions of uniqueness and existence of weak solution). But as it is well known, the maximum principle even for the simple case of hyperbolic equation (one dimensional wave equation \cite{Protter}) are quite different from those for elliptic and parabolic cases,  for which it is a natural fact, such a way a role of characteristics curves and surfaces becomes evident in  the situation of hyperbolic type PDEs.

 We call the angle of characteristics slop the solution to the equation $-\tan\varphi_j=\lambda_j$, and angle between $j-$ and $k-$ characteristics: $\varphi_k-\varphi_j\neq \pi  l,\, l\in Z,$     where $\lambda_j\neq\pm i$ are real and prime roots of the characteristics equation, $j,\,k=1,\,2,\,3,\,4$.

 Most of these equations serve as  mathematical models of many physical processes and attract the interest of researchers. The most famous of them are elasticity beam equations (Timoshenko beam equations with and without internal damping) \cite{Capsoni}, short laser pulse equation \cite{Fichera}, equations which describe the structures are subjected to moving loads, and equation of Euler-Bernoulli beam resting on two-parameter Pasternak foundation and subjected to a moving load or mass \cite{Uzz} and others.

Due to evident practice application, these models need more exact tools for studying, and as consequence, to attract fundamental knowledge. As usual, most of these models are studied by the analytical-numerical methods (Galerkin's methods).

On the other hand, the maximum principle is an efficient tool for fundamental knowledge of PDEs: the range of problems, which maximum principle allows to study, belongs to a class of quite actual problems of well-posedness of so-called general boundary-value problems for higher-order differential equations originating from the works by L. H\"ormander and M.Vishik who used the theory of extensions to prove the existence of well-posed boundary-value problems for linear differential equations of arbitrary order with constant complex coefficients in a bounded domain with smooth boundary. This theory got its present-day development in the works by G. Grubb \cite{Grubb}, L.Hormander \cite{H1}, and A. Posilicano \cite{Pos}. Later, the problem of well-posedness of boundary-value problems for various types of second order differential equations was studied by V. Burskii and A. Zhedanov \cite{Bursk1}, \cite{Bursk2} which developed a method of traces associated with a differential operator and applied this method to establish the Poncelet, Abel and Goursat problems, and by I. Kmit \cite{Kmit}. In the previous works of author (see \cite{Bur3}) there have been developed qualitative methods of studying  Cauchy problems and nonstandard in the case of hyperbolic equations Dirichlet and Neumann problems for the linear fourth-order equations (moreover, for an equation  of any even order $2m,\, m\geq 2,$ ) with the help of operator methods (L-traces, theory of extension, moment problem, method of duality equation-domain and others), \cite{Bur2}.

As concern maximum principle, at the present time there are not any results for the fourth order equations even in linear case. As it was mentioned above, the maximum principle even for the simple case of one dimensional wave equation \cite{Protter}, and for the second-order telegraph equation \cite{Maw1}---\cite{Ortega} are quite different from those for elliptic and parabolic cases.
In the monograph of Protter and  Weinberger \cite{Protter} there was shown that solutions of hyperbolic equations and inequalities do not exhibit the classical formulation of maximum principle. Even in the simplest case of the wave equation in two independent variables $u_{tt}-u_{xx}=0$  the maximum of a nonconstant solution $u = \sin x\sin t$ in a rectangle domain $\{(x,t) :\, x \in [0, \pi],\, t \in [0, \pi] \}$  occurs at an interior point $ \left(\frac{\pi}{2},\,\frac{\pi}{2}\right)$. In Chapter 4 \cite{Protter} the maximum principle for linear second hyperbolic equations of general type, with variable coefficients has also been obtained for Cauchy problems and boundary value problems on characteristics (Goursat problem).

Following R.Ortega, A.Robles-Perez \cite{Ortega}, we introduce the definition of a "weak form" of the maximum principle, which is used for the hyperbolic equations, which will be used later.
\vskip5mm
 {\it Definition 1.} \cite{Ortega}
 Let $ L= Lu$ be linear differential operator, acting on functions $u:\,D\to R$, in some domain $D$. These functions will belong to the certain family $ B,$ which includes boundary conditions or others requirements. It is said that $ L$ satisfies the maximum principle, if
$$ L\geq 0,\,\,u\in B,$$
implies $u\geq 0$ in $D$.
 \vskip5mm

 In further works of these authors (see \cite{Maw1}, \cite{Maw2},\cite{Maw3}) there was studied  the maximum principle for weak bounded twice periodical solutions from the space $L^{\infty}$ of the telegraph equation with parameter $\lambda$ in lower term, one-, two-, and -tree dimensional spaces, and which includes the cases of variables coefficients. The precise condition for $\lambda$ under which the maximum principle still valid was font.  There was also introduced a method of upper and lower solutions  associated with the nonlinear equation, which allows to obtain the analogous results (uniqueness, existence and regularity theorems) for the  telegraph equations with external nonlinear forcing, applying maximum principle. There was considered also the case when the external forcing belongs to a certain space of measures.

 The maximum principle for general quasilinear hyperbolic systems with dissipation was proved by Kong De Xing \cite{Kong}. There were given two estimates for the solution to the general quasilinear hyperbolic system and introduced the concept of dissipation (strong dissipation and weak dissipation), then state some maximum principles of quasilinear hyperbolic systems with dissipation. Using the maximum principle there were reproved the existence and uniqueness theorems of the global smooth solution to the Cauchy problem for considered quasilinear hyperbolic system.

 So, the problem to prove the maximum principle for the weak solutions stills more complicated and at that time becomes more interesting in the case of fourth-order hyperbolic equations, especially, in the case of non-classical boundary value problems with weak-regularity data.  There are no results on maximum principle even for model case of linear 2-dimensional fourth-order hyperbolic equations with the constant coefficients and without lower terms. Moreover, we can not use the term of usual traces in the cases of initial data of weak regularity, and we come to the notions of $L-$traces, the traces, which associated with differential operator. Let us remind (see, for example, \cite {Bursk1}), that $L-$ traces exist for the weak solutions from space $L^2$ even in the situations when classical notions of traces does not work for such solutions.

 \section{Statement of the problem and auxiliary definitions.}
Let us start to establish the maximum principle for the weak solutions to the Cauchy problem
for the Eq.(\ref{eq:01}) in some admissible planar domain. It is  expected, that in hyperbolic case the  characteristics of the equations play a crucial role.

Let $C_j,\,j=1,\,2,\,3,\,4$ be characteristics,  $\Gamma_0:=\{x_1\in [a,\,b],\,x_2=0\},$ and  define domain $\Omega$ as a domain which is restricted by the characteristics $C_j,\,j=1,\,2,\,3,\,4$ and $\Gamma_0.$ Following \cite{Protter} we will call below the domain $\Omega$ as characteristics domain. For the second order hyperbolic equations (see \cite{Protter}), $\Omega$ is characteristics triangle, in the case of fourth-order equations with constant coefficients, that is existence of 4 different and real characteristics lines $C_j,\,j=1,\,2,\,3,\,4$, $\Omega$ is a characteristic pentagon.

Consider also the following Cauchy problem for the Eq. (\ref{eq:01}) on $\Gamma_0:$
\begin{equation}\label{eq:02}
u|_{\Gamma_0}=\varphi(x),\,u^{\prime}_{\nu}|_{\Gamma_0}=\psi(x),\,u^{\prime\prime}_{\nu\nu}|_{\Gamma_0}=\sigma(x),\,u^{\prime\prime\prime}_{\nu\nu\nu}|_{\Gamma_0}=
\chi(x),
\end{equation}
where $\varphi,\,\psi,\,\sigma$ and $\chi$ are given weak regular functions on $\Gamma_0,$ in general case $ \varphi,\,\psi,\,\sigma,\,\chi\in L^2(\Gamma_0),\,\nu-$ is outer normal of $\Gamma_0.$

{\it Definition 2.} We call a domain $D:=\{(x_1,\,x_2):\,x_1\in (-\infty,\,+\infty),\,x_2>0\}$ in the half-plane $x_2 > 0$ an admissible domain if
it has the property that for each point $C\in D$ the corresponding characteristics domain $\Omega$ is also in $D$.  More generally, $D$ is admissible, if it is the finite or countable union of characteristics domains.

We choose some arbitrary point $C\in D$ in admissible plane domain $D,$  draw through this point two arbitrary characteristics, $C_1$ and $C_2$. Another two characteristics ($C_3$ and $C_4$) we draw through the ends $a$ and $b$ of initial line $\Gamma_0$. We determine the points $O_1$ and $O_2$ as intersections of $C_1,\,C_3$ and $C_2,\,C_4$ correspondingly: $O_1=C_1\cap C_3,\,O_2=C_2\cap C_4.$ Such a way, domain $\Omega$ is a pentagon $aO_1CO_2b.$

Establishment of the maximum principle in this situation allows us to obtain a local properties of the solution to Cauchy problem (\ref{eq:01})---(\ref{eq:02}) on the arbitrary interior point $C\in D.$

We will consider the weak solution to the problem (\ref{eq:01})---(\ref{eq:02}) from the $D(L)$, domain of definition of maximal operator, associated with the differential operation $L$ in Eq.(\ref{eq:01}). Following \cite{Bur3},\,\cite{Grubb} and \cite{H1}, we remind the corresponding definitions.

In the bounded domain $\Omega$ we consider the linear differential operation ${\mathcal L}$ of the order $m,\,m\geq 2,$ and formally adjoint ${\mathcal L}^+$:
\begin{equation}\label{eq:03}
{\mathcal L}(D_x)=\sum\limits_{|\alpha|\leq m}a_{\alpha}D^{\alpha},\,{\mathcal L}^+(D_x)=\sum\limits_{|\alpha|\leq m}D^{\alpha}(a_{\alpha}),
\end{equation}
where $\alpha=(\alpha_1,\,\alpha_2,...\alpha_n),\,|\alpha|=\alpha_1+\alpha_2+...+\alpha_n$ is multi-index. Note, that for Eq. (\ref{eq:01}) $n=2,\,m=4.$

{\it Definition 3. Minimum operator.} \cite{Bur3}.
Let us consider the differential operation (\ref{eq:03}) on functions from the space $C_0^{\infty}(\Omega).$ The minimum operator $L_0$ is called the extension of the operation from $C_0^{\infty}(\Omega)$ to the set $D(L_0):=\overline{C_0^{\infty}(\Omega)}.$ The closure is realized in the norm of the graph of operator $L$: $||u||_L^2:=||u||_{L_2(\Omega)}^2+||Lu||^2_{L_2(\Omega)}.$

{\it Definition 4. Maximum operator.} \cite{Bur3}.
The maximum operator $L$ is define as the restriction of the differential operation $  {\mathcal L}(D_x)$ to the set $D(L):=\{u\in L^2(\Omega):\,Lu\in L^2(\Omega)\}.$

{\it Definition 5.} \cite{Bur3}.
The operator $\tilde L$ is define as the extension of the minimum operator $L_0,$ to the set $D(\tilde L):=\overline{C^{\infty}(\bar\Omega)}.$

{\it Definition 6. Regular operator.} \cite{Bur3}.
The maximum operator is called regular if  $D(L)=D(\tilde L).$

For example, for elliptic case $D(\tilde L)=H^4(\Omega),$ $D(L_0)=\stackrel{0}{H^4}(\Omega),$ the  Hilbert Sobolev space of fourthly weak differentiable functions from $L^2(\Omega)$ (\cite{Bur5}).

Analogously, can be introduced the operators $L^+,\tilde L^+,$ and $L_0^+,$ associated with the formally adjoint operation ${\mathcal L}^+$.

The definition of a weak solution to the problem (\ref{eq:01})---(\ref{eq:02}) from the space $D(L)$ is closely connected with the notion of $L-$ traces, that is traces, which are associated with the differential operator $L$.

{\it Definition 7. L-traces.} \cite{Bur5}.
Assume, that for a function $u\in D(\tilde L)$ there exist linear continuous functionals $L_{(p)}u$ over the space $H^{m-p-1/2}(\partial\Omega),\,p=0,1,2...,m-1$, such that the following equality is satisfied:
\begin{equation}\label{eq:04}
(Lu,v)_{L^2(\Omega)}-(u,L^+v)_{L^2(\Omega)}=\sum\limits_{j=0}^{m-1}(L_{(m-1-j)}u,\,\partial^{(j)}_{\nu}v),
\end{equation}
for any functions $v\in H^m(\Omega).$
The functionals $L_{(p)}u$ is called the $L_{(p)}-$ traces of the function $u\in D(\tilde L).$ Here $(\cdot,\,\cdot)_{L^2(\Omega)}$ is a scalar product in Hilbert space $L^2(\Omega)$.

For $L^2-$ solutions the notion of $L_{(p)}-$ traces can be realized by the following way.

{\it Definition 8. } The distributions $L_{(p)}u\in H^{-p-\frac{1}{1}}(\partial\Omega),\,p=0,...,\,m-1,$ are called the $p-$th $L-$traces of the function $u\in D(L)$ on $\partial\Omega$, if the following identity is true
\begin{equation}\label{eq:04-1}
\int\limits_{\Omega}\left(Lu\cdot\overline{v}-u\cdot\overline{L^+v}\right)\,dx=\sum\limits_{j=0}^{m-1}<L_{(m-1-j)}u,\,\partial^{(j)}_{\nu}v>_{\partial\Omega}.
\end{equation}
for any functions $v\in H^m(\Omega).$

For example, relations for $L-$traces of the solution $u\in D(L)$ to the fourth-order equation $Lu=f\in L^2(\Omega)$ have the form:
$$\sum\limits_{j=0}^{3}<L_{(3-j)}u,\,\partial^{(j)}_{\nu}v>_{\partial\Omega}=\int\limits_{\Omega}f\cdot\overline{v}\,dx,$$
for all $v\in Ker\,L^+\cap H^m(\Omega).$

Finally, we are going to the definition of the weak solution to the problem (\ref{eq:01})---(\ref{eq:02}):

{\it Definition 9.} We will call the function $u\in D(L)$ a weak solution to the Cauchy problem (\ref{eq:01})--(\ref{eq:02}), if it satisfies to the following integral identity
\begin{equation}\label{eq:05}
(f,\,v)_{L^2(\Omega)}-(u,\,L^+v)_{L^2(\Omega)}=\sum\limits_{j=0}^{3}<L_{(3-j)}u,\,\partial^{(j)}_{\nu}v>_{\partial\Omega},
\end{equation}
for any functions $v\in H^m(\Omega).$
The functionals $L_{(p)}u$ is called the $L_{(p)}-$ traces of the function $u$, $p=0,\,1,\,2,\,3,$ and completely determined by the initial functions $\varphi,\,\psi,\,\sigma,\,\chi$ by the following way:
\begin{equation*}
L_{(0)}u=-L(x)u|_{\partial\Omega}=-L(\nu)\varphi;
\end{equation*}
\begin{equation*}
L_{(1)}u=L(\nu)\psi+\alpha_1\varphi_{\tau}^{\prime}+\alpha_2\varphi;
\end{equation*}
\begin{equation}\label{eq:06}
L_{(2)}u=-L(\nu)\sigma+\beta_1\psi_{\tau}^{\prime}+\beta_2\psi+\beta_3\varphi_{\tau\tau}^{\prime\prime}+\beta_4\varphi_{\tau}^{\prime}+\beta_5\varphi;
\end{equation}
\begin{equation*}
L_{(3)}u=L(\nu)\chi+\delta_1\varphi_{\tau\tau\tau}^{\prime\prime\prime}+\delta_2\sigma+\delta_3\psi_{\tau\tau}^{\prime\prime}+\delta_4\psi_{\tau}^{\prime}+
\delta_5\psi+\delta_6\varphi_{\tau\tau}^{\prime\prime}+\delta_7\varphi_{\tau}^{\prime}+\delta_8\varphi.
\end{equation*}
Here $\alpha_i,\,i=1,\,2,\,\beta_j,\,j=1,\,2,...,\,5,$ and $\delta_k,\,k=1,\,...,\,9$ are smooth functions, completely determined by the coefficients of the Eq.(\ref{eq:01}).

{\it Remark 1.} We can use a general form of the operators $\gamma_j$ in the left-hand side of the identity (\ref{eq:05}) instead of operators of differentiation $\partial^{(j)}_{\nu}v$. Indeed, we define $\gamma_j=p_j\gamma,$ where $\gamma: \,u\in H^m(\Omega)\to (u|_{\partial\Omega},\,...,\,u_{\nu}^{(m-1)}|_{\partial\Omega})\in H^{(m)}=H^{m-1/2}(\partial\Omega)\times H^{m-3/2}(\partial\Omega)\times...\times H^{1/2}(\partial\Omega),$ and $p_j:\,H^{(m)}\to H^{m-j-1/2}(\partial\Omega)-$ projection.

As it has been mentioned above,  some examples show (see \cite{Bursk1}) that in the general case the solutions $u\in D(L)$ do
not exist ordinary traces in the sense of distributions even for the simplest hyperbolic equations.
 Indeed, for the wave equation $Lu =\frac{\partial^2 u}{\partial x_1\partial x_2} = 0$ in the unit disk $K:\, |x|=1$, the solution $u(x) = (1-x_1^2)^{-\frac{5}{2}}$
belongs to $L^2(K)$, but $<u|_{\partial K}, 1>_{\partial K} = \infty$ it means that $\lim_{r\to 1-0}\int\limits_{|x|=r}u(x)ds_x=\infty$, such a way
the trace $u|_{\partial K}$ does not exist even as a distribution. However, for every solution $u\in L^2(K)$  the $L_{(0)}-$trace $L_{(0)}u := -L(x) u(x)|_{|x|=1}=-x_1x_2 u(x)|_{|x|=1}\in L^2(\partial K).$  Likewise, $L_{(1)}-$ trace $L_{(1)}u$ exists for every $u\in L^2(K)$:
\begin{equation*}
L_{(1)}u=\left(L(x)u^{\prime}_{\nu}+L^{\prime}_{\tau}u^{\prime}_{\tau}+\frac{1}{2}L^{\prime\prime}_{\tau\tau}u\right)|_{\partial K}\in H^{-\frac{3}{2}}(\partial K).
\end{equation*}
where $\tau$ is the angular coordinate and $u^{\prime}_{\tau}$
is the tangential derivative, and $L(x)=x_1x_2-$ symbol of the operator $L=\frac{\partial^2 }{\partial x_1\partial x_2}$.

\section{Maximum principle for the weak solutions of Cauchy problem.}
We prove here the first simple case: the maximum principle for the weak solution of the Cauchy problem (\ref{eq:01})---(\ref{eq:02}) in admissible plane domain $\Omega,$ restricted by the different and  non-congruent characteristics $C_j,\,j=1,\,2,...,\,4$ and initial line $\Gamma_0$.

{\it Theorem 1. Maximum principle.} Let $u\in D(L)$ satisfy the following inequalities:

\begin{equation}\label{eq:07}
 Lu=f\leq 0,\,\,\,x\in D,
 \end{equation}
and
 \begin{equation}
L_{(0)}u\mid_{\Gamma_0}\geq 0,\,L_{(1)}u|_{\Gamma_0}\geq 0,\,L_{(2)}u|_{\Gamma_0}\geq 0,\,
L_{(3)}u|_{\Gamma_0}\geq 0,\label{eq:08}
\end{equation}
then
$u\leq 0$ in $D.$
\vskip3.5mm
{\it Proof.}
1. First of all prove the statement for smooth solutions $u\in C^{\infty}(\overline{\Omega}).$

Due to homogeneity of the symbol in Eq. (\ref{eq:01}), $L(\xi)=a_0\xi_1^4+a_1\xi_1^3\xi_2+a_2\xi_1^2\xi_2^2+a_3\xi_1\xi_2^3+a_4\xi_2^4=\\ <\xi,\,a^1><\xi,\,a^2><\xi,\,a^3><\xi,\,a^4>,\,\xi=(\xi_1,\,\xi_2)\in  R^2,$ we can rewrite this equation in the following form:
\begin{equation}
<\nabla,\,a^1><\nabla,\,a^2><\nabla,\,a^3><\nabla,\,a^4>u=f(x).\label{eq:09}
\end{equation}
The vectors $a^j=(a_1^j,\,a_2^j),\,j=1,\,2,\,3,\,4$ are determined by the coefficients $a_i,\,i=0,\,1,\,2,\,3,\,4, $ and $<a,\,b>=a_1\bar b_1+a_2\bar b_2$ is a scalar product in $ C^2$. It is easy to see that vector $a^j$ is a tangent vector of $j-$th characteristic, slope $\varphi_j$ of which is determined by $-\tan\varphi_j=\lambda_j,\,j=1,\,2,\,3,\,4.$ In what follows, we also consider the vectors $\tilde a^j=(-\bar a_2^j,\,\bar a_1^j),\,j=1,\,2,\,3,\,4.$ It is obvious that $<\tilde a^j,\,a^j>=0,$ so $\tilde a^j$ is a normal vector of $j-$th characteristic.

Use the definitions 7 and 9 for the case  $m=4,$ that is fourth-order operator in Eq. (\ref{eq:01}), and domain  $\Omega, $ which is restricted by the characteristics $C_j,\,j=1,\,2,\,3,\,4$ and $\Gamma_0:$
\begin{equation*}
\int\limits_{\Omega}\{Lu\cdot \bar v-u\cdot\overline{L^+v}\}dx=\sum\limits_{k=0}^{3}\int\limits_{\partial\Omega}L_{(3-k)}u\cdot \partial^{(k)}_{\nu}v\,ds=
\end{equation*}
\begin{equation*}
=\sum\limits_{k=0}^{3}\int\limits_{C_1}L_{(3-k)}u\cdot \partial^{(k)}_{\nu}v\,ds+\sum\limits_{k=0}^{3}\int\limits_{C_2}L_{(3-k)}u\cdot \partial^{(k)}_{\nu}v\,ds+\sum\limits_{k=0}^{3}\int\limits_{C_3}L_{(3-k)}u\cdot \partial^{(k)}_{\nu}v\,ds+\sum\limits_{k=0}^{3}\int\limits_{C_4}L_{(3-k)}u\cdot \partial^{(k)}_{\nu}v\,ds+
\end{equation*}
\begin{equation}
+\sum\limits_{k=0}^{3}\int\limits_{\Gamma_0}L_{(3-k)}u\cdot \partial^{(k)}_{\nu}v\,ds.\label{eq:10}
\end{equation}
Using the representation (\ref{eq:09}), we arrive to
$$\int\limits_{\Omega}Lu\cdot \bar v\,dx=\int\limits_{\Omega}<\nabla,\,a^1><\nabla,\,a^2><\nabla,\,a^3><\nabla,\,a^4>u\cdot\bar v\,dx=$$
$$\int\limits_{\partial\Omega}<\nu,\,a^1>\cdot<\nabla,\,a^2><\nabla,\,a^3><\nabla,\,a^4>u\cdot\bar v\,ds-$$
$$\int\limits_{\Omega}<\nabla,\,a^2><\nabla,\,a^3><\nabla,\,a^4>u\cdot\overline{<\nabla,\,a^1>v}\,dx.$$
Integrating by parts further, we obtain:
$$\int\limits_{\Omega}Lu\cdot \bar v\,dx=
\int\limits_{\partial\Omega}<\nu,\,a^1>\cdot<\nabla,\,a^2><\nabla,\,a^3><\nabla,\,a^4>u\cdot\bar v\,ds-$$
$$\int\limits_{\partial\Omega}<\nu,\,a^2>\cdot<\nabla,\,a^3><\nabla,\,a^4>u\cdot\overline{<\nabla,\,a^1>v}\,ds+$$
$$+\int\limits_{\partial\Omega}<\nu,\,a^3>\cdot<\nabla,\,a^4>u\cdot\overline{<\nabla,\,a^2><\nabla,\,a^1>v}\,ds-$$
$$-\int\limits_{\partial\Omega}<\nu,\,a^4>\cdot u\cdot\overline{<\nabla,\,a^3><\nabla,\,a^2><\nabla,\,a^1>v}\,ds+$$
$$+\int\limits_{\Omega}u\cdot\overline{<\nabla,\,a^4><\nabla,\,a^3><\nabla,\,a^2><\nabla,\,a^1>v}\,dx.$$
Since $<\nabla,\,a^4><\nabla,\,a^3><\nabla,\,a^2><\nabla,\,a^1>v=L^+v,$ and determining
$$\tilde L_{(0)}u:=<\nu,\,a^4>u,\,\,\tilde L_{(1)}u:=<\nu,\,a^3>\cdot<\nabla,\,a^4>u,$$
$$\tilde L_{(2)}u:=<\nu,\,a^2>\cdot<\nabla,\,a^3><\nabla,\,a^4>u,$$
 $$\tilde L_{(3)}u=L_{(3)}u=<\nu,\,a^1>\cdot<\nabla,\,a^2><\nabla,\,a^3><\nabla,\,a^4>u$$ that are
analogues of $L-$traces from the formula (\ref{eq:10}). Such a way we have
$$\int\limits_{\Omega}\{Lu\cdot \bar v-u\cdot\overline{L^+v}\}\,dx=
\int\limits_{\partial\Omega}L_{(3)}u\cdot\bar v\,ds-\int\limits_{\partial\Omega}\tilde L_{(2)}u\cdot\overline{<\nabla,\,a^1>v}\,ds+$$
\begin{equation}\label{eq:11}
+\int\limits_{\partial\Omega}\tilde L_{(1)}u\cdot\overline{<\nabla,\,a^2><\nabla,\,a^1>v}\,ds
-\int\limits_{\partial\Omega}\tilde L_{(0)} u\cdot\overline{<\nabla,\,a^3><\nabla,\,a^2><\nabla,\,a^1>v}\,ds.
\end{equation}
Difference between formulas (\ref{eq:10}) and (\ref{eq:11}) is that natural $L_{(3-k)}$ traces in (\ref{eq:10}) are multiplied by $k-$ derivative by outer normal $\nu$ of truncated function $v:\,\partial^{(k)}_{\nu}v,$ on the other hand, in (\ref{eq:11}) we determined by $\tilde L_{(3-k)}$ some expressions which multiplied by differential operators $L_{k}^+v$ of order $k$ and which can serve as  analogous of natural  $L_{(3-k)}$ traces, $k=0,\,1,\,2,\,3.$ So, in the (\ref{eq:11})
$$L^+_1v:=<\nabla,\,a^1>v,\,L^+_2v:=<\nabla,\,a^2><\nabla,\,a^1>v,$$
$$L_0^+v=v,\,L^+_3v:=<\nabla,\,a^3><\nabla,\,a^2><\nabla,\,a^1>v.$$
Let  $v\in Ker L^+$ in (\ref{eq:11}) and calculate $L-$ traces on $\partial\Omega=C_1\cup C_2\cup C_3\cup C_4\cup\Gamma_0$. For instance, for $L_{(3)}u$ we obtain: $L_{(3)}u=<\nu,\,a^1><\nabla,\,a^2><\nabla,\,a^3><\nabla,\,a^4>u, $ and use that $<\nabla,\,a^j>u=$\\$<\nu,\,a^j>u^{\prime}_{\nu}+<\tau,\,a^j>u^{\prime}_{\tau},\,j=1,\,2,\,3,\,4,$ where $\nu-$ normal vector, $\tau-$ tangent vector. Due to presence the product $<\nu,\,a^1>,$  $L_{(3)}u=0$ on characteristic $C_1,$ normal vector $\tilde a^1$ of which is orthogonal to the vector $a^1$. On the other parts of $\partial\Omega$ there will vanish the terms containing $<\nu,a^j>$ on $C_j$.   After that\
$$\int\limits_{\partial\Omega}<\nu,\,a^1><\nabla,\,a^2><\nabla,\,a^3><\nabla,\,a^4>u=\int\limits_{\Gamma_0}L_{(3)}u\,ds+$$
$$<\tilde a^2,\,a^1><a^2,\,a^2><\tilde a^2,\,a^3><\tilde a^2,\,a^4>\int\limits_{C_2}u_{\nu\nu\tau}\,ds+$$
$$<\tilde a^3,\,a^1><\tilde a^3,\,a^2>< a^3,\,a^3><\tilde a^3,\,a^4>\int\limits_{C_3}u_{\nu\nu\tau}\,ds+$$
$$<\tilde a^4,\,a^1><\tilde a^4,\,a^2><\tilde a^4,\,a^3>< a^4,\,a^4>\int\limits_{C_4}u_{\nu\nu\tau}\,ds+$$
$$\left\{<\tilde a^2,\,a^1><a^2,\,a^2><\tilde a^2,\,a^3><a^2,\, a^4>+<\tilde a^2,\,a^1><a^2,\,a^2>< a^2,\,a^3><\tilde a^2,\, a^4>\right\}\int\limits_{C_2}u_{\tau\tau\nu}\,ds+$$
$$ \left\{<\tilde a^3,\,a^1><\tilde a^3,\,a^2><a^3,\,a^3><a^3,\,a^4>+<\tilde a^3,\,a^1><a^3,\,a^2>< a^3,\,a^3><\tilde a^3,\, a^4>\right\}\int\limits_{C_3}u_{\tau\tau\nu}\,ds+$$
$$\left\{<\tilde a^4,\,a^1><\tilde a^4,\,a^2><a^4,\,a^3><a^4,\, a^4>+<\tilde a^4,\,a^1><a^4,\,a^2>< \tilde a^4,\,a^3><a^4,\, a^4>\right\}\int\limits_{C_4}u_{\tau\tau\nu}\,ds+$$
$$<\tilde a^2,\,a^1><a^2,\,a^2>< a^2,\,a^3><a^2,\,a^4>\int\limits_{C_2}u_{\tau\tau\tau}\,ds+$$
$$<\tilde a^3,\,a^1><a^3,\,a^2>< a^3,\,a^3><a^3,\,a^4>\int\limits_{C_3}u_{\tau\tau\tau}\,ds+$$
$$<\tilde a^4,\,a^1><a^4,\,a^2>< a^4,\,a^3><a^4,\,a^4>\int\limits_{C_4}u_{\tau\tau\tau}\,ds+
\alpha_{4,1}\int\limits_{C_2}u_{\nu\nu}\,ds+\alpha_{4,2}\int\limits_{C_3}u_{\nu\nu}\,ds+\alpha_{4,3}\int\limits_{C_4}u_{\nu\nu}\,ds+$$
$$\alpha_{5,1}\int\limits_{C_2}u_{\nu\tau}\,ds+\alpha_{5,2}\int\limits_{C_3}u_{\nu\tau}\,ds+\alpha_{5,3}\int\limits_{C_4}u_{\nu\tau}\,ds
+\alpha_{6,1}\int\limits_{C_2}u_{\tau\tau}\,ds+\alpha_{6,2}\int\limits_{C_3}u_{\tau\tau}\,ds+\alpha_{6,3}\int\limits_{C_4}u_{\tau\tau}\,ds+$$
$$\alpha_{7,1}\int\limits_{C_2}u_{\nu}\,ds+\alpha_{7,2}\int\limits_{C_3}u_{\nu}\,ds+\alpha_{7,3}\int\limits_{C_4}u_{\nu}\,ds+
\alpha_{8,1}\int\limits_{C_2}u_{\tau}\,ds+\alpha_{8,2}\int\limits_{C_3}u_{\tau}\,ds+\alpha_{8,3}\int\limits_{C_4}u_{\tau}\,ds.$$
Here correspondent coefficients $\alpha_{i,j}$ were numerated as follows: first index $i$ indicates the derivative of $u$: \\ $1)\,u_{\nu\nu\tau},\,2)\,u_{\nu\tau\tau},\,3)\,u_{\tau\tau\tau},\,4)\, u_{\nu\nu},\,5)\,u_{\nu\tau},\,6)\,u_{\tau\tau},\,7)\, u_{\nu}\, 8)\,u_{\tau},$ the second index $j$ indicates $j+1-$th characteristic, $j=1,\,2,\,3.$
So, now the  formula (\ref{eq:10}) has the form:
$$\int\limits_{\Omega}Lu\,dx=\int\limits_{\Gamma_0}L_{(3)}u\,ds+
\alpha_{1,1}\int\limits_{C_2}u_{\nu\nu\tau}\,ds+\alpha_{1,2}\int\limits_{C_3}u_{\nu\nu\tau}\,ds+\alpha_{1,3}\int\limits_{C_4}u_{\nu\nu\tau}\,ds+$$
$$\alpha_{2,1}\int\limits_{C_2}u_{\tau\tau\nu}\,ds+\alpha_{2,2}\int\limits_{C_3}u_{\tau\tau\nu}\,ds+\alpha_{2,3}\int\limits_{C_4}u_{\tau\tau\nu}\,ds+$$
$$\alpha_{3,1}\int\limits_{C_2}u_{\tau\tau\tau}\,ds+\alpha_{3,2}\int\limits_{C_3}u_{\tau\tau\tau}\,ds+\alpha_{3,3}\int\limits_{C_4}u_{\tau\tau\tau}\,ds+$$
$$\alpha_{4,1}\int\limits_{C_2}u_{\nu\nu}\,ds+\alpha_{4,2}\int\limits_{C_3}u_{\nu\nu}\,ds+\alpha_{4,3}\int\limits_{C_4}u_{\nu\nu}\,ds+$$
$$\alpha_{5,1}\int\limits_{C_2}u_{\nu\tau}\,ds+\alpha_{5,2}\int\limits_{C_3}u_{\nu\tau}\,ds+\alpha_{5,3}\int\limits_{C_4}u_{\nu\tau}\,ds+
\alpha_{6,1}\int\limits_{C_2}u_{\tau\tau}\,ds+\alpha_{6,2}\int\limits_{C_3}u_{\tau\tau}\,ds+\alpha_{6,3}\int\limits_{C_4}u_{\tau\tau}\,ds+$$
$$\alpha_{7,1}\int\limits_{C_2}u_{\nu}\,ds+\alpha_{7,2}\int\limits_{C_3}u_{\nu}\,ds+\alpha_{7,3}\int\limits_{C_4}u_{\nu}\,ds+
\alpha_{8,1}\int\limits_{C_2}u_{\tau}\,ds+\alpha_{8,2}\int\limits_{C_3}u_{\tau}\,ds+\alpha_{8,3}\int\limits_{C_4}u_{\tau}\,ds.$$
Coefficients $\alpha_{i,j}$ are constant and depend on only from Eq. (\ref{eq:01}) coefficients $a_0,\,a_1,\,a_2,\,a_3,\,a_4.$ By analogous way we calculate others $L-$ traces, $L_{(0)}u,\,L_{(1)}u$ and $L_{(2)}u.$

The value of the function $u$ at the point $C\in D,\, u(C)$ we estimate from the last equality, integrating by the characteristics $C_1$ and $C_2$ and using conditions (\ref{eq:02}), (\ref{eq:06})---(\ref{eq:08}). Since, the chosen point $C\in D$ is arbitrary, we arrive at $u\leq 0$ in $D.$

2. For solutions $u\in D(L)$ the statement of the Theorem 1 follows from the conditions:
$$ \overline{C^{\infty}(\overline{\Omega}})=D(L),$$
and
$$ \overline{C^{\infty}(\overline{\Omega}})=D(L^+).$$
These conditions hold for operators with constant coefficients in domain convex with respect characteristics (see \cite{H1}).

The Theorem 1 is proved.

{\it Remark 2.} The weak form of the maximum principle for  $u\in L^2(\Omega)$ can be derived not only for the solutions of the Cauchy problem  (\ref{eq:02}) but for all linear differential operator problems $Lu=F\in L^2(\Omega),$  with condition $\overline{{\rm Im}\, L^+}=L^2(\Omega)$ and with the constant coefficients.

Indeed, using conditions (\ref{eq:07}), (\ref{eq:08}), and the definition 9 we obtain
$$\int\limits_{\Omega}u\cdot\overline{L^+v}\,dx\leq 0,$$
for all $v\in H^m(\Omega).$ If $\overline{{\rm Im}\, L^+}=L^2(\Omega),$ then we have
$$\int\limits_{\Omega}u\cdot\overline{w}\,dx\leq 0,$$
for any $w\in L^2(\Omega).$
The last inequality can serve as a weak maximum principle for solutions to the boundary value problems from $L^2(\Omega).$

{\it Remark 3.} In the case of classical solution of the Cauchy problem for the second order hyperbolic equations of the general form with the constant coefficients the statement of the Theorem 1 coincides with the result of \cite{Protter}.  In this case conditions (\ref{eq:08}) have usual form without using the notion of $L-$traces (see \cite{Protter}):
\begin{equation*}
u|_{\Gamma_0}\leq 0,\,\,u^{\prime}_{\nu}|_{\Gamma_0}\leq 0.
\end{equation*}

\section*{Acknowledgments.}
The author thanks to Prof. Iryna Kmit for the drawing the author’s attention to the problem of the maximum principle for high-order hyperbolic PDEs and systems,that allowed to apply developed earlier methods of investigations of the fourth-order PDEs. This work is supported by the Volkswagen Foundation (the project numbers are A131968 and 9C624) and by the Ministry of Education and Science of Ukraine (project number is 0121U109525).



%
%

\bigskip

\end{document}